\documentclass{amsart}
\usepackage{amssymb}
\usepackage{amsfonts}
\usepackage{amsmath}
\usepackage[legalpaper,bookmarks=true,colorlinks=true,linkcolor=blue,citecolor=blue]{hyperref}
\usepackage{graphicx}%
\usepackage{fancyhdr}
\usepackage{color}
\usepackage[mathlines]{lineno}
\usepackage{lscape}
\usepackage{epsfig}
\usepackage{natbib}
\usepackage{geometry}
\usepackage{tgbonum}
\usepackage[]{listings}
\fontfamily{qcr}\selectfont

\newtheorem{theorem}{Theorem}
\theoremstyle{plain}

\newtheorem{proposition}{Proposition}

\numberwithin{equation}{section}

\newcommand{\Bin}{\bigskip \noindent}

\newcommand{\Ni}{\noindent}

\begin{document}
\Large
\title[$\ell^{\infty}$ Poisson invariance principles from two Poisson limit theorems and extension]{
$\ell^{\infty}$ Poisson invariance principles from two classical Poisson limit theorems and extension to non-stationary independent sequences}

\author{Gane Samb Lo} \author{Aladji Babacar Niang}\author{Amadou Ball} 

\begin{abstract}  The simple Lévy Poisson process and scaled forms are explicitly constructed from partial sums of independent and identically distributed random variables and from sums of non-stationary independent random variables. For the latter, the weak limits are scaled Poisson processes. The method proposed here prepares generalizations to dependent data, to associated data in the first place.\\

\noindent $^{\dag}$ Gane Samb Lo.\\
LERSTAD, Gaston Berger University, Saint-Louis, S\'en\'egal (main affiliation).\newline
LSTA, Pierre and Marie Curie University, Paris VI, France.\newline
AUST - African University of Sciences and Technology, Abuja, Nigeria\\
Imhotep Mathematical Center (IMC), imhotepsciences.org\\
gane-samb.lo@edu.ugb.sn, gslo@aust.edu.ng, ganesamblo@ganesamblo.net\\
Permanent address : 1178 Evanston Dr NW T3P 0J9,Calgary, Alberta, Canada.\\

\noindent $^{\dag\dag}$ Aladji Babacar Niang\\
LERSTAD, Gaston Berger University, Saint-Louis, S\'en\'egal.\\
Email: niang.aladji-babacar@ugb.edu.sn, aladjibacar93@gmail.com\\
Imhotep Mathematical Center (IMC), imhotepsciences.org\\

\noindent $^{\dag\dag\dag}$ Amadou Ball\\
LERSTAD, Gaston Berger University, Saint-Louis, S\'en\'egal (main affiliation).\newline
fany.ngom@ugb.edu.sn\\

\noindent $^{\dag}$ Gane Samb Lo.\\
LERSTAD, Gaston Berger University, Saint-Louis, S\'en\'egal (main affiliation).\newline
LSTA, Pierre and Marie Curie University, Paris VI, France.\newline
AUST - African University of Sciences and Technology, Abuja, Nigeria\\
Imhotep Mathematical Center (IMC), imhotepsciences.org\\
gane-samb.lo@edu.ugb.sn, gslo@aust.edu.ng, ganesamblo@ganesamblo.net\\
Permanent address : 1178 Evanston Dr NW T3P 0J9,Calgary, Alberta, Canada.\\

\noindent \textbf{Keywords}: central limit theorem for arrays of random variables; infinitely divisible laws; simple and scaled Poisson L\'evy processes; weak convergence of stochastic processes in the space of bounded functions on $(0,1)$.\\

\Ni\textbf{AMS 2010 Mathematics Subject Classification}: 60F17; 60E07.
\end{abstract}
\maketitle

\newpage

\noindent \textbf{R\'{e}sum\'{e}.} Les processus poisonniens de Lévy simples ou re-échelonnés sont construits explicitely à partir de sommes partielles de de suites de variables aléatoires indépendantes et identiquement distribuées et de suites non-stationaires indépendentes. Pour ces dernières, les lois limites faibles sont des processus de Poisson composés. Cette étude prépare le terrain à des généralisations aux données dépendantes, aux variables associées dans un premier temps.\\

\Ni \textbf{The authors}.\\

\Ni \textbf{Aladji babacar Niang}, M.Sc., is preparing a Ph.D. dissertation under the supervision of the fourth authors at Gaster Berger University, SENEGAL.\\

\Ni \textbf{Gane Samb Lo}, Ph.D. and Doctorate in Sciences, is full professor at Gaston Berger of Saint-Louis (SENEGAL), at the African University of Science and Technology (AUST), NIGERIA. He is affiliated to LSTA, Pierre et Marie Curie University, FRANCE. He is the founder of Imhotep Mathematical Center\\

\Ni \textbf{Ch\'erif Mactar Mamadou Traor\'e}, M.Sc., is preparing a Ph.D. dissertation under the supervision of the second authors at Gaster Berger University, SENEGAL.\\

\Ni \textbf{Amadou Ball}, M.Sc., has prepared a M>Sc degree under the supervision of the second authors at Gaston Berger University, SENEGAL. Most of the work of this paper was done in his dissertation.\\

\section{Introduction}

\Ni In this paper, we provide invariance principles, also called functional limit theorems, for Poisson weak limits for two classical results in probability theory but also for their recent generalizations. Let us begin by describing the central limit theorems for which functional laws have to be established.

\subsection{Poisson weak limits} \label{sec_01_ss_02}
\Ni The approximation of a sequence of binomial probability laws $\left(\mathcal{B}(n,p_{n})\right)_{n\geq 1}$ associated to a sequence of random variables $(Z_n)_{n\geq 1}$ [such that the sequence of probabilities $(p_n)_{n\geq 1}$ converges to zero and $np_n \rightarrow \lambda >0$ as $n\rightarrow +\infty$] to a Poisson law $\mathcal{P}(\lambda)$ is a classical and easy-to-prove result in probability theory. This approximation is very important in some real-life situations, especially in lack of powerful computers. In this simple case, the distribution of each $Z_n$ is a convolution product of $n$ independent and identically distributed (\textit{iid}) Bernoulli laws $\mathcal{B}(p_n)$-associated to the random variables
$\left\{(X_{j,n})_{1\leq j \leq n}, \ n\geq 1\right\}$, \textit{i.e.}, $Z_n=X_{1,n}+X_{2,n}+\cdots+X_{n,n}$. A parallel theory also exists for sums of corrected geometric laws. When we depart from the \textit{iid} assumption, the problem may get more and rapidly, even if the independence assumption is still required. The situation becomes more interesting if the random variables $X_{j,n}$ are non-stationary and dependent. Generalizations of these results have been given recently by \cite{apwl-nnLo-ins}. It happens that asymptotic laws of sums of random variables are closely related to invariance principles or functional weak limits which in turn may lead to L\'evy stochastic processes which are so important in many areas of applications, in mathematical finance for example.\\

\Bin In that view, the aforementioned Asymptotic Poisson Weak Limits (\textit{APWL}) should be associated with Invariance principles (\textit{IP}) or Functional Limit Theorems (\textit{FLT}) leading to scaled L\'evy process, both for the classical cases and for the extended non-stationary approach. Among L\'evy processes, two are iconic: Brownian motions and Poisson processes.\\

\Bin So our aim is to establish \textit{FLT}'s for the extended Poisson weak limit laws in \cite{apwl-nnLo-ins} and by the way to re-establish \textit{FLT}'s for the two classical cases, since the statements of the later cases can be hardly found in the most common literature. \\ 

\Ni To have more details on our objectives, we recall the \textit{APWL} under consideration in Subsection \ref{sec_01_ss_02}. Next we make a quick introduction
to \textit{FLT}'s in the space $\ell^{\infty}(0,1)$ in Subsection \ref{sec_01_ss_03} with the main tools to be used there. \\

\subsection{The Asymptotic Poisson Weak Laws} \label{sec_01_ss_02}
\Ni We recall the two results for stationary Bernoulli and corrected geometric sequences of random variables.

\begin{proposition} \label{piidBinomial} Let $(X_{n})_{n\geq 0}$ be a sequence of random variables in some probability space $(\Omega,\mathcal{A}, \mathbb{P})$ such that:\\

\Ni 1) $\forall n\geq 1 $, $X_n \sim \mathcal{B}(n,p_{n})$,  $n\geq 1 $;\\

\Ni 2) $p_{n} \rightarrow 0$ and $np_{n} \rightarrow \lambda \in \mathbb{R}_{+}\setminus\{0\}$ as $n\rightarrow + \infty$.\\

\Ni Then 

$$
X_{n} \rightsquigarrow  \mathcal{P}(\lambda).
$$
\end{proposition}

\Bin Next, we have: \\

\begin{proposition}\label{niidBinomial}
Let $(X_{n})_{n\geq 0}$ be a sequence of random variables in some probability space $(\Omega,\mathcal{A}, \mathbb{P})$ such that: \\

\Ni 1) $\forall n\geq 1 $, $X_n \sim \mathcal{N}\mathcal{B}(n,p_{n})$,  $n\geq 1 $; \\

\Ni 2) $(1-p_{n}) \rightarrow 0 $ and $ n(1-p_{n}) \rightarrow \lambda \in \mathbb{R}_{+}\setminus \{0\} $ as $n\rightarrow + \infty$.\\

\Ni Then 

$$
X_{n}-n \rightsquigarrow  \mathcal{P}(\lambda).
$$
\end{proposition}

\Bin We also recall the two generalizations of the above mentioned results for non-stationary Bernoulli and corrected geometric sequences of random variables (see \cite{apwl-nnLo-ins}).

\begin{theorem}\label{pbinomialINS} Let $X=\biggr\{ \{X_{k,n}, \ 1 \leq k \leq k_n=k(n)\}, \ n\geq 1\biggr\}$ be an array by-row-independent Bernoulli random variables, that is:\\

\Ni (1) $\forall n\geq 1$, $\forall 1\leq k \leq k(n)$, $X_{k,n}\sim \mathcal{B}(p_{k,n})$, with $0<p_{k,n}<1$ and:\\

\Ni (2) $\sup_{1\leq k \leq k(n)} p_{k,n} \rightarrow 0$;\\

\Ni (3) $\sum_{1\leq k \leq k(n)} p_{k,n} \rightarrow \lambda \in ]0, \ +\infty[$.\\

\Ni Then we have

$$
S_n[X]:=\sum_{k=1}^{k(n)} X_{k,n} \rightsquigarrow \mathcal{P}(\lambda).
$$

\end{theorem}

\Bin Next: \\

\begin{theorem}\label{nbinomialINS}  Let $X=\biggr\{ \{X_{k,n}, \ 1 \leq k \leq k_n=k(n)\}, \ n\geq 1\biggr\}$ be an array by-row-independent corrected Geometric  random variables, that is: \\

\Ni (1) $\forall n\geq 1$, $\forall 1\leq k \leq k(n)$, $X_{k,n}\sim \mathcal{G}^{\ast}(p_{k,n})$, with $0<p_{k,n}=1-q_{k,n}<1$ and:\\

\Ni (2) $\sup_{1\leq k \leq k(n)} q_{k,n} \rightarrow 0$;\\

\Ni (3) $\sum_{1\leq k \leq k(n)} q_{k,n} \rightarrow \lambda \in ]0, \ +\infty[$.\\

\Ni Then we have

$$
S_n[X]:=\sum_{k=1}^{k(n)} X_{k,n} \rightsquigarrow \mathcal{P}(\lambda).
$$
\end{theorem}

\Bin Now, let us do a quick introduction on simple functional laws.\\

\subsection{A quick introduction to the \textit{FLT} in $\ell^{\infty}(0,1)$} \label{sec_01_ss_03}
\noindent The general principle of invariance principles is the following. Let $X_{1}$, $X_{2}$,$\cdots$ be a sequence of real-valued centered random variables, defined on the same probability space, with finite variances, that is $\sigma_i^2=\mathbb{E}\left\vert X_{i}\right\vert ^{2}<\infty$. For each $n\geq 1$, set

$$
s_n^2=\sigma_1^2+\cdots + \sigma_n^2
$$

\bigskip\noindent and

\begin{equation*}
S_{n}=X_{1}+\cdots+X_{n}.
\end{equation*}

\bigskip\noindent For $0\leq t\leq 1$ and $n\geq 1$, put
  
\begin{equation*}
Y_{n}(t)=\frac{S_{\left[ nt\right] }}{s_n},
\end{equation*}

\bigskip\noindent where, for any real $u$, $[u]$ stands for the integer part of $u$, which is the greatest integer less or equal to $u$.\\

\noindent In the \textit{iid} case with variance one, $s_n^2=n$ and the sequence $\{Y_{n}(t), \ 0\leq t \leq 1\}$ is studied in the space $D(0,1)$ of functions of the first kind endowed with the Skorohod metric. It may be proved that

\begin{equation}
\{Y_{n}(t), \ 0\leq t \leq 1\} \rightsquigarrow \{W(t), \ 0\leq t \leq 1\}, \ as \ n\rightarrow +\infty, \label{FLT}
\end{equation}

\bigskip\noindent meaning that $\{Y_{n}(t), \ 0\leq t \leq 1\}$ weakly converges to a Wiener process (Brownian Motion) $\{W(t), \ 0\leq t \leq 1\}$ in the sense of the Skorohod topology.\\

\noindent This result is the starting point of a considerable research trying to have extensions of that result, labeled as \textit{functional laws} or  \textit{invariance principles}.\\

\noindent \textbf{In each part of this document}, the notation given in this introduction will be used and eventually adapted.\\

\Ni For a long time, the weak law \eqref{FLT} proved in the space $\mathcal{C}(0,1)$ of continuous functions on (0,1) or in the space $\mathcal{D}(0,1)$ of functions defined on (0,1) having left and right limits at each point and having at most a countable number of discontinuity points. The spaces $\mathcal{C}(0,1)$ and $\mathcal{D}(0,1)$, when endowed with the supremum norm and the Skorohod metric are complete and separable spaces (Polish spaces). Most importantly, suprema of stochastic processes in them are measurable. However, handling such weak limits in them, specially in $\mathcal{D}(0,1)$ is very hard. To find a less complicated way, the space $\ell^{\infty}(T)$ of bounded functions equipped with the supremum norm is used where $T$ is some set, and here we use $T=[a,b]$ and frequently $T=[0,\ 1]$. To avoid the measurability problem, exterior and interior integrals with respect to a probability measure, and exterior and interior probabilities are introduced. A very advanced round up of weak convergence in $\ell^{\infty}(0,1)$ is available in \cite{vaart}.\\

\Ni Now, we can state our precise objective, which is to prove that for all four theorems above, we will have \textit{FLT}'s in the form

\begin{equation}
\{Y_{n}(t), \ 0\leq t \leq 1\} \rightsquigarrow \{N(\lambda a(t)), \ 0\leq t \leq 1\}, \ in \ \ell^{\infty}(0,1) \ as \ n\rightarrow +\infty, \label{FLT2}
\end{equation}

\Bin where $N(\lambda a(\circ))$ is scaled Poisson process of intensity $\lambda$ and that $a(t)\equiv t$ for the first two theorems.\\

\Ni Finally, before we give our results, let us describe the method of finding \textit{FTL}'s in $\ell^{\infty}(0,1)$.\\

\subsection{General tools for weak laws in $\ell^{\infty}(0,1)$}
\Ni The invariance principle results are proved by first establishing the finite-dimensional convergence to a given stochastic process and second, by showing that the sequence of stochastic process is asymptotically tight. We will use the following theorem. \\

\begin{theorem} \label{workingTool_01} Let $a$ and $b$ be two real numbers such that $a<b$ and $T=[a,b]$. For the sequence $(Y_{n})_{n\geq 1} \subset \ell^{\infty} \left([a,b]\right)$ converges to a tight stochastic process $W \in \ell^{\infty}\left([a,b]\right)$, it is sufficient that :\\

\noindent (a) the finite dimensional margins of  $(Y_{n})_{n\geq 1}$ converge to those of $W$\\

\noindent  and\\

\noindent (b) $(Y_{n})_{n\geq 1}$ is asymptotically tight.
\end{theorem}

\bigskip \noindent The asymptotically tightness is characterized as follows: Given that each margin $(Y_{n}(t))_{n\geq 1}$, $t \in T$, is asymptotically tight, the stochastic process $(Y_{n})_{n\geq 1}$ is tight if and only if there exists a semi-metric $\rho$ on $T=[a,b]$ such that $(T,\rho)$ is totally bounded and for all $\eta>0$,

\begin{equation}
\lim_{\delta \downarrow 0} \limsup_{n\rightarrow +\infty} \mathbb{P}^\ast\left(\sup_{\rho(t,s)<\delta} \left|Y_n(t)-Y_n(s)\right|\geq \eta\right)=0. \ \ \label{ip_tightness_charac}
\end{equation}

\bigskip\noindent (See \cite{vaart} or \cite{ips-wcib-ang}, Chapter 3). \\

\Ni However in most cases where $T$ is a bounded interval of $\mathbb{R}$, we try to have that tightness in the simple case where $\rho(s,t)=|t-s|$. \\

\section{\textit{FLT} related \textit{APWL}} \label{sec_02}

\Ni Let us consider the stationary case first. We have the two following theorems.\\

\begin{theorem} \label{FLT-piidBinomial} Let $(X_{n})_{n\geq 0}$ be a sequence of random variables in some probability space $(\Omega,\mathcal{A}, \mathbb{P})$ such that:\\

\Ni 1) $\forall n\geq 1 $, $X_n=X_{1,n} + \cdots + X_{n,n}$, with each $X_{i,n} \sim \mathcal{B}(p_{n})$, \ $1\leq i\leq n$;\\

\Ni 2) $p_{n} \rightarrow 0$ and $np_{n} \rightarrow \lambda \in \mathbb{R}_{+}\setminus\{0\}$ as $n\rightarrow + \infty$.\\

\Ni Let us consider the sequence of stochastic processes

$$
\biggr\{ \left\{Y_n(t), \ 0\leq t \leq 1\right\}, \  n\geq 1 \biggr\}=\biggr\{ \left\{X_{[nt]}, \ 0\leq t \leq 1\right\}, \  n\geq 1 \biggr\}
$$

\Bin with $X_{[nt]}=0$ for $0\leq nt <1$. Then we have

\begin{equation}
Y_{n}(\circ) \rightsquigarrow N(\circ) \ in \  \ell^{\infty}(0,1),  \label{FLT2SB}
\end{equation}

\Bin where $N(\circ)$ is a Poisson process of intensity $\lambda$.
\end{theorem}

\Bin We also have: \\

\begin{theorem}\label{FLT-niidBinomial}
Let $(X_{n})_{n\geq 0}$ be a sequence of random variables in some probability space $(\Omega,\mathcal{A}, \mathbb{P})$ such that: \\

\Ni 1) $\forall n\geq 1 $, $X_n=X_{1,n} + \cdots + X_{n,n}$, with each $X_{i,n} \sim \mathcal{G}(p_{n})$, \ $1\leq i\leq n$;\\

\Ni 2) $(1-p_{n}) \rightarrow 0 $ and $ n(1-p_{n}) \rightarrow \lambda \in \mathbb{R}_{+}\setminus \{0\} $ as $n\rightarrow + \infty$.\\

\Ni Let us denote $Z_n=X_n-n$, $n\geq 1$ and consider the sequence of stochastic processes

$$
\biggr\{ \left\{Y_n(t), \ 0\leq t \leq 1\right\}, \  n\geq 1 \biggr\}=\biggr\{ \left\{Z_{[nt]}, \ 0\leq t \leq 1\right\}, \  n\geq 1 \biggr\},
$$

\Bin with $Z_{[nt]}=0$ for $0\leq nt <1$. Then we have

\begin{equation}
Y_{n}(\circ) \rightsquigarrow N(\circ) \ in \  \ell^{\infty}(0,1),  \label{FLT2SG}
\end{equation}

\Bin where $N(\circ)$ is a Poisson process of intensity $\lambda$.
\end{theorem}

\Bin \textbf{Proofs of theorems}.\\

\Ni \textbf{Proof of Theorem \ref{FLT-piidBinomial}}. We are going to apply Theorem \ref{workingTool_01}. Let us proceed to two steps.\\

\Ni \textbf{Step 1}. Let us begin by the finite distribution convergences. Since, we have for any $0\leq t \leq 1$,

$$
Y_n(t)=X_{[nt]}=X_{1,n}+\cdots+ X_{[nt],n} \sim \mathcal{B}([nt],p_n), \ for \ t\geq 1/n,
$$

\Bin we have, for $n\geq 1$ big enough,

$$
[nt]p_n=(np_n) ([nt]/n) \to \lambda t
$$

\Bin and by Proposition \ref{piidBinomial}, $Y_n(t) \rightsquigarrow \mathcal{P}(\lambda t)$. Now let $k\geq 0$ and $0=t_0<t_1<\cdots<t_k$. Let $n\geq 1$ large enough (say $n\geq n_0$) to ensure that all $nt_j\geq 1$, $1\leq j \leq k$ and $[nt_1]<\cdots<[nt_k]$. Let us define

$$
Z_n=\left(Y_n(t_1), \ Y_n(t_2)-Y_n(t_1),\ldots,Y_n(t_k)-Y_n(t_{k-1})\right), \ n\geq n_0.
$$

\Bin It is clear that:\\

\Ni (1) $Z_n$ has independent components; \\

\Ni (2) For each $1\leq j\leq k$, we have

$$
Z_n(t_j)=Y_n(t_j)-Y_n(t_{j-1})=X_{[nt_{j-1}]+1}+\cdots+X_{[nt_{j}]},
$$

\Bin and the number of terms in $Z_n(t_j)$ is  $m_n(j)=[nt_{j}]-[nt_{j-1}]$ and satisfies $m_n(j)p_n \rightarrow \lambda(t_j-t_{j-1})$. Hence, by Proposition \ref{piidBinomial},
 
$$
Z_n(j) \rightsquigarrow \mathcal{P}\left(\lambda (t_j-t_{j-1})\right).\\
$$

\Bin So, $Z_n=(Z_{n}(t_1), \cdots, Z_{n}(t_k))$ has independent components which weakly converge to marginal laws. By the general Slutsky rule, $Z_n$ weakly converges to the product law of those marginals laws, that is 

$$
Z_n \rightsquigarrow (Z_{t_1}, \cdots, Z_{t_k}),
$$

\Bin where the $Z_{t_1}$, $\cdots$, $Z_{t_k}$ are independent and $Z_{t_j} \sim \mathcal{P}(\lambda (t_{j}-t_{j-1}))$. So, by decrementing, we obtain that
$(Y_n(t_1),\cdots,Y_n(t_k))$ converges to $(N(t_1,\lambda), \cdots, N(t_k,\lambda))$, where $N(\circ,\lambda)$ is a simple Poisson process. $\square$\\

\Ni \textbf{Step 2}. Let $\delta>0$ and let, for $n\geq 1$

$$
A_n(\delta)= \sup_{(s,t)\in [0,1]^2, \ |s-t|<\delta} |Y_n(s)-Y_n(t)|.
$$

\Bin We take $\delta \in ]0,1[$ and so $m(n,\delta):=\left[n(\delta+1/n)\right]$ is small against $n$. Here we will exploit that any $Y_n(t)$ is a sum of constant-signed random variables. We have the following facts.\\

\Ni (a) For any $(s,t)\in [0,1]^2$ such that $|s-t|<\delta$, $|Y_n(s)-Y_n(t)|$ is of the form $X_{j,n}+\cdots+X_{j+h,n}$, where $h\leq m(n,\delta)-1$.\\

\Ni (b) So for each $j \in \{1,\cdots,n-m(n,\delta)+1\}$, fixed, $X_{j,n}+\cdots+X_{j+h,n}$, with $h\leq m(n,\delta)-1$, is bounded by 

$$
Z_{j,n}:=X_{j,n}+\cdots+X_{j+m(n,\delta)-1},n. 
$$

\Bin Let us denote $N=n-m(n,\delta)+1$. For $j\leq N$, the quantity  $X_{j,n}+\cdots+X_{j+h,n}$, with $h\leq m(n,\delta)-1$, is still bounded by $Z_N:=\sup_{1\leq j \leq N} Z_{j,n}$.\\

\Ni (c) Finally, we have 

$$
A_n(\delta)\leq \sup_{1\leq j \leq N}  Z_{j,n}.
$$

\Ni We set $B_h=(\sup_{1\leq j \leq N}  Z_{j,n} = Z_{h,n})$, $1\leq h \leq N$. It is clear that

$$
\Omega = \bigcup_{1\leq h \leq N} B_h=\sum_{1\leq h \leq N} B^{\prime}_h,
$$

\Bin with $B^{\prime}_1=B_1$, $B^{\prime}_2=B_1^cB_2$, $B^{\prime}_h=B_1^c \cdots B_{h-1}^c B_h$ for $h\geq 1$. We get, for any $\eta>0$,

\begin{eqnarray*}
\mathbb{P}(A_n >\eta)&\leq& \mathbb{P}(\{\sup_{1\leq j \leq N} Z_{j,n}\}>\eta)\\
&=&\sum_{h=1}^{N} \mathbb{P}\left(\left\{ \sup_{1\leq j \leq N} Z_{j,n} >\eta\right\} \cap B^{\prime}_h\right)\\
&=& \sum_{h=1}^{N} \mathbb{P}\left(\left\{ \sup_{1\leq j \leq N} Z_{j,n} >\eta\right\}  / B^{\prime}_h\right) \ \mathbb{P}( B^{\prime}_h)\\
&=& \sum_{h=1}^{N} \mathbb{P}(\{Z_{h,n} >\eta\}) \ \mathbb{P}( B^{\prime}_h)\\
&\leq & \frac{1}{\eta} \sum_{h=1}^{N} \mathbb{E}(Z_{h,n}) \ \mathbb{P}( B^{\prime}_h)\\
&=& \frac{m(n,\delta)p_n}{\eta}  \sum_{h=1}^{N} \mathbb{P}( B^{\prime}_h)\\
&=& \frac{m(n,\delta)p_n}{\eta}= \frac{\delta (\lambda+o(1))}{\eta}.
\end{eqnarray*}

\Bin We conclude that, for any $\eta>0$,

$$
\lim_{\delta\downarrow 0} \limsup_{n\rightarrow +\infty} \mathbb{P}\biggr(\biggr(\sup_{(s,t)\in [0,1]^2, \ |s-t|<\delta} |Y_n(s)-Y_n(t)|\biggr) > \eta\biggr)=0.
$$

\Bin $\square$\\

\Bin
\textbf{Proof of Theorem \ref{FLT-niidBinomial}}. That proof closely follows that of Theorem \ref{FLT-piidBinomial}. There is no need to give the details. 
$\blacksquare$

\section{\textit{FLT} for non-stationary Bernoulli or corrected geometrical laws} \label{sec_03}

\Bin To have clues on how to extend the above results to a non-stationary scheme, we may draw some facts from hypotheses of Theorem \ref{FLT-piidBinomial} for example. There, we add $k_n(t)=[nt]$, $0\leq t\leq 1$, $n\geq 1$.  For $0\leq s <t$,

$$
Y_n(t)-Y_n(s)=X_{k_n(s)+1,n}+\cdots + X_{k_n(t),n} \sim \mathcal{B}(k_n(t)-k_n(s),p_n).
$$

\Bin For $a(t)=t$, \ $0\leq t \leq 1$, by setting $\Delta a(s,t)=a(t)-a(s)$,

$$
\mathbb{E}(Y_n(t)-Y_n(s))=(k_n(t)-k_n(s))p_n=: \Delta p_n(s,t) \rightarrow \lambda \Delta a(s,t)=\lambda \Delta(s,t), 
$$

\Bin uniformly in $(s,t) \in [0,1]^2$. This simple analysis suggests generalizing the above \textit{FLT} in the following way: require a sequence of functions 
$(k_n(t))_{0\leq t \leq 1}$, $n\geq 1$, which are non-decreasing with $k_n(0)=0$ and $k_n(1)=n$ and a uniformly right-continuous and non-decreasing function $a(t)$ of $t\in [0,1]$, with $a(0)=1-a(1)=0$ such that for $0\leq s <t$,

$$
\lim_{n\rightarrow +\infty} \sup_{(s,t) \in [0,1]^2} \left|\biggr(p_{k_n(s)+1,n} + \cdots + p_{k_n(t),n}\biggr)- \lambda \Delta a(s,t)\right|=0
$$

\Bin and

$$
\lim_{\delta \rightarrow 0} \sup_{(s,t) \in [0,1]^2, |s-t|<\delta} \Delta a(s,t)=0.
$$

\Bin Exploiting these ideas leads to the two \textit{FLT}'s for non stationary data.\\

\subsection{Non-stationary sequences of Bernoulli random variables} \label{sec_02_ss_01}

\begin{theorem} \label{FLT-pBinomialNS} Let $X=\biggr\{ \{X_{k,n}, \ 1 \leq k \leq k_n=k(n)\}, \ n\geq 1\biggr\}$ (with $k_n\rightarrow +\infty$ as $n\rightarrow +\infty$) be an array by-row-independent Bernoulli random variables, that is: \\

\Ni (1) $\forall n\geq 1$, $\forall 1\leq k \leq k(n)$, $X_{k,n}\sim \mathcal{B}(p_{k,n})$, with $0<p_{k,n}<1$.\\

\Bin Suppose that we have the following other assumptions.\\

\Ni (2) $\sup_{1\leq k \leq k(n)} p_{k,n} \rightarrow 0$;\\

\Ni (3) There exist:\\

\Bin (3a) a sequence of non-decreasing functions $(k_n(t))_{0\leq t \leq 1}$, such that $k_n(0)=0$ and $k_n(1)=k_n$, $n\geq 1$, (satisfying : for any 
$0\leq s <t$, we have $k_n(s)<k_n(t)$ for $n$ large enough) \\

\Ni and: \\

\Ni (3b) a uniformly right-continuous and non-decreasing function $a(t)$ of $t\in [0,1]$ with $a(0)=1-a(1)=0$, \textit{i.e}, 

$$
\lim_{\delta \rightarrow 0} \sup_{t \in [0,1[} \Delta a(t,t+\delta)=0; \ \ (F3a)
$$

\Bin with $\Delta a(s,t)=a(t)-a(s)$ for $0\leq s <t$;  \\

\Ni Moreover, for

$$
\Delta p_n(s,t)=p_{k_n(s)+1,n} + \cdots + p_{k_n(t),n}, 
$$

\Bin we have

$$
\limsup_{n\rightarrow +\infty} \left| \Delta p_n(s,t) - \lambda \Delta a(s,t)\right|=0 \ \ (F3b)
$$

\Bin and, for $\delta>0$ (small enough) and for $m(n,\delta)=k_n(t+\delta)-k_n(t)$,

$$
\limsup_{n\rightarrow +\infty} \sup_{1\leq j \leq k_n-m(n,\delta)+1} \left|p_{j,n}+\cdots+p_{j+m(n,\delta)-1,n} - \lambda \delta\right|=0. \ \ (F3c)
$$

\Bin Then the stochastic process

$$
\biggr\{ \left\{Y_n(t), \ 0\leq t \leq 1\right\}, \  n\geq 1 \biggr\}=\biggr\{ \left\{X_{1,n}+\cdots+X_{k_n(t),n}, \ 0\leq t \leq 1\right\}, \  n\geq 1 \biggr\}
$$

\Bin weakly converges in $\ell^{\infty}(0,1)$ to  the compound Poisson process $N(a(\circ),\lambda)$ of intensity $\lambda>0$.

\end{theorem}

\Bin \textbf{Proof of Theorem of \ref{FLT-pBinomialNS}}. Let us proceed by steps.\\

\Ni \textbf{Step 1}. Finite dimensional convergence. Let $r>1$ and $0=t_0<t_1<\cdots<t_r\leq 1$. Let for $1\leq j \leq r$,

$$
Z_{j,n}=Y_n(t_j)-Y_n(t_{j-1})=X_{k_n(t_{j-1})+1,n}+\cdots+X_{k_n(t_{j}),n}.
$$ 

\Bin Since all the assumptions of Theorem \ref{pbinomialINS} hold for each $Z_{j,n}$, we have

$$
\forall j \in \{1,\cdots,r\}, \ Z_{j,n} \rightsquigarrow \mathcal{P}(\lambda \Delta a(t_{j-1},t_j)). 
$$

\Bin So, $Z_n=(Z_{1,n}, \cdots, Z_{r,n})$ has independent components which weakly converge to marginal laws. By the general Slutsky rule, $Z_n$ weakly converges to the product law of those marginals laws, that is, 

$$
Z_n \rightsquigarrow (Z_1, \cdots, Z_r),
$$

\Bin where the $Z_1$, $\cdots$, $Z_r$ are independent and $Z_j \sim \mathcal{P}(\lambda \Delta a(t_{j-1},t_j))$. So, by decrementing, we obtain that
$(Y_n(t_1),\cdots,Y_n(t_r))$ converges to the $(N(a(t_1),\lambda), \cdots, N(a(t_r),\lambda))$, where $N(a(\circ),\lambda)$ is a scaled Poisson process by $(a(\circ))$.\\

\Ni \textbf{Step 2}. We closely follow the proof of Theorem \ref{FLT-piidBinomial} and re-use its notations. Let $\delta>0$ and let, for $n\geq 1$

$$
A_n(\delta)= \sup_{(s,t)\in [0,1]^2, \ |s-t|<\delta} |Y_n(s)-Y_n(t)|.
$$

\Bin We have:\\

\Ni (a) For any $(s,t)\in [0,1]^2$, such that $|s-t|<\delta$, $|Y_n(s)-Y_n(t)|$ is of the form $X_{j,n}+\cdots+X_{j+h,n}$, where 

$$
h\leq m(n,\delta)-1, \ \ with \ \ m(n,\delta)=k_n(t+\delta)-k_n(\delta).
$$

\Bin (b) So for each $j \in \{1,\cdots, k_n-m(n,\delta)+1\}$, fixed, $X_{j,n}+\cdots+X_{j+h,n}$, with $h\leq m(n,\delta)-1$, is bounded by 

$$
Z_{j,n}=:X_{j,n}+\cdots+X_{j+m(n,\delta)-1},n. 
$$

\Bin Let us denote $N=k_n-m(n,\delta)+1$. For $j\leq N$, the quantity  $X_{j,n}+\cdots+X_{j+h,n}$, with $h\leq m(n,\delta)-1$, is still bounded by $Z_N=:\sup_{1\leq j \leq N}  Z_{j,n}$.\\

\Ni (c) Finally, we have 

$$
A_n(\delta)\leq \sup_{1\leq j \leq N}  Z_{j,n}.
$$

\Bin We set $B_h=(\sup_{1\leq j \leq N}  Z_j = Z_h)$, $1\leq h \leq N$. Surely, we have

$$
\Omega = \bigcup_{1\leq h \leq N} B_h=\sum_{1\leq h \leq N} B^{\prime}_h
$$

\Bin with $B^{\prime}_1=B_1$, $B^{\prime}_2=B_1^cB_2$, $B^{\prime}_h=B_1^c \cdots B_{h-1}^c B_h$ for $h\geq 1$. We get, for any $\eta>0$,

\begin{eqnarray*}
\mathbb{P}(A_n(\delta) >\eta)&\leq& \mathbb{P}(\{\sup_{1\leq j \leq N} Z_{j,n}\}>\eta)\\
&=&\sum_{h=1}^{N} \mathbb{P}\left(\left\{ \sup_{1\leq j \leq N} Z_{j,n} >\eta\right\} \cap B^{\prime}_h\right)\\
&=& \sum_{h=1}^{N} \mathbb{P}\left(\left\{ \sup_{1\leq j \leq N} Z_{j,n} >\eta\right\}  / B^{\prime}_h\right) \ \mathbb{P}( B^{\prime}_h)\\
&=& \sum_{h=1}^{N} \mathbb{P}(\{Z_{h,n} >\eta\}) \ \mathbb{P}( B^{\prime}_h)\\
&\leq & \frac{1}{\eta} \sum_{h=1}^{N} \mathbb{E}(Z_{h,n}) \ \mathbb{P}( B^{\prime}_h)\\
&=& \frac{1}{\eta}  \sum_{h=1}^{N} \biggr(p_{h,n}+\cdots+p_{h+m(n,\delta)-1,n}\biggr) \ \mathbb{P}( B^{\prime}_h) \ \ (L26)\\
&=& \frac{\lambda \Delta(t,t+\delta)}{\eta}  \sum_{h=1}^{N} \ \mathbb{P}( B^{\prime}_h)\\
&=& \frac{\lambda \Delta(t,t+\delta)}{\eta},
\end{eqnarray*}

\Bin where we applied Assumption (F3c) in Line (L26). By applying Assumption (F3a), we get, for any $\eta>0$,

$$
\lim_{\delta\downarrow 0} \limsup_{n\rightarrow +\infty} \mathbb{P}\biggr(\biggr(\sup_{(s,t)\in [0,1]^2, \ |s-t|<\delta} |Y_n(s)-Y_n(t)|\biggr) > \eta\biggr)=0.
$$

\Bin $\square$\\

\subsection{Non-stationary sequences of Corrected geometric variables} \label{sec_02_ss_02}

\Ni The following \textit{FLT} is also valid for sums of corrected geometric random variables. The proof will be omitted because it is very similar of the proof of Theorem \ref{FLT-pBinomialNS}. \\

\begin{theorem} \label{FLT-nBinomialNS} Let $X=\biggr\{ \{X_{k,n}, \ 1 \leq k \leq k_n=k(n)\}, \ n\geq 1\biggr\}$ (with $k_n\rightarrow +\infty$ as $n\rightarrow +\infty$) be an array by-row-independent corrected geometric random variables, that is: \\

\Ni (1) $\forall n\geq 1$, $\forall 1\leq k \leq k(n)$, $X_{k,n}\sim \mathcal{G}^{\ast}(p_{k,n})$, with $0<p_{k,n}<1$.\\

\Ni Suppose that we have the following other assumptions.\\

\Ni (2) $\sup_{1\leq k \leq k(n)} q_{k,n} \rightarrow 0$.\\

\Ni (3) There exist:\\

\Bin (3a) a sequence of non-decreasing functions $(k_n(t))_{0\leq t \leq 1}$ such that $k_n(0)=0$ and $k_n(1)=k_n$, $n\geq 1$ (satisfying : for any 
$0\leq s <t$, we have $k_n(s)<k_n(t)$ for $n$ large enough) \\

\Ni and: \\

\Ni (3b) a uniformly right-continuous and non-decreasing functions $a(t)$ of $t\in [0,1]$ with $a(0)=1-a(1)=0$, \textit{i.e.},

$$
\lim_{\delta \rightarrow 0} \sup_{t \in [0,1]} \Delta a(t,t+\delta)=0 \ \ (G3a)
$$

\Bin with $\Delta a(s,t)=a(t)-a(s)$ for $0\leq s <t$. Moreover, by setting  

$$
\Delta q_n(s,t)=q_{k_n(s)+1,n} + \cdots + q_{k_n(t),n}, 
$$

\Bin we have

$$
\limsup_{n\rightarrow +\infty} \left| \Delta q_n(s,t) - \lambda \Delta a(s,t)\right|=0 \ \ (G3b)
$$

\Bin and, for $\delta>0$ (small enough) and for $m(n,\delta)=k_n(t+\delta)-k_n(t)$,

$$
\limsup_{n\rightarrow +\infty} \sup_{1\leq j \leq k_n-m(n,\delta)+1} \left|q_{j,n}+\cdots+q_{j+m(n,\delta)-1,n} - \lambda \delta\right|=0. \ \ (G3c)
$$

\Bin Then the sequence of stochastic processes

$$
\biggr\{ \left\{Y_n(t), \ 0\leq t \leq 1\right\}, \  n\geq 1 \biggr\}=\biggr\{ \left\{X_{1,n}+\cdots+X_{k_n(t),n}, \ 0\leq t \leq 1\right\}, \  n\geq 1 \biggr\}
$$

\Bin weakly converges in $\ell^{\infty}(0,1)$ to the compound Poisson process $N(a(\circ),\lambda)$ of intensity $\lambda>0$. \\

\end{theorem}
 
%======================== -1

\subsection{Examples} \label{sec-example}

\Ni Let us provide some examples for which the assumptions of Theorems \ref{FLT-pBinomialNS} and \ref{FLT-nBinomialNS}, specifically (F3b), (F3c), (G3b) and (G3c) hold. Actually, we illustrate conditions of the first theorem only.\\

\Ni \textbf{1. Stationary case}. We already said that (F3b) and (F3c) hold in the stationary case with $k_n=n$, $k_n(t)=[tk_n]$, $p_{k,n}=p_n$ and $a(t)=t$. This still hold even if $k_n$ is an arbitrary sequence converging to $+\infty$.\\

\Ni \textbf{2. Checking (F3c)}. It seems that (F3c) is more difficult to get. The example used to realize (F3b) in the next example seems not to make hold (F3c).\\

\Ni Actually, (F3c) broadly means the partial sums $p_{j,n}+\cdots+p_{j+\ell,n}$, for $\ell>0$, are stationary, \textit{i.e.}, nearly stationary. They would be exactly stationary if they were all equal to $p_{1,n}+\cdots+p_{1+\ell,n}$. In other words, (F3c) means that the partial sums $p_{j,n}+\cdots+p_{j+\ell,n}$ asymptotically behave like $p_{1,n}+\cdots+p_{1+\ell,n}$, independently of $j$, as $n\rightarrow +\infty$, meaning that they are asymptotically stationary.\\

\Ni Let us give some examples of such sequences.\\

\Ni \textbf{3. Non trivial example for (F3b) holds}. Let us consider an array $\{\{X_{k,n, \ 1\leq k\leq \ell(n)}\}, \ n\geq 1\}$ with $\ell(n)\rightarrow +\infty$. Let $\varepsilon \in ]0,1[$. Let $(k_n)_{n\geq 1}$ be a sequence converging to $+\infty$ such that $[k_n^{1+\varepsilon}]\leq \ell(n)$. Let us define

$$
k_n(t)=\left[k_n^{\varepsilon+t}\right], \ n\geq 1, \ t\in [0,1],
$$

\Bin for $\gamma$ standing for the Euler number,

$$
b_n=\sum_{1\leq k \leq k_n} k^{-1}=\log k_n + \gamma + o(1).
$$

\Bin Now, suppose that for any $1<k\leq K_n=\left[ k_n^{1+\varepsilon}\right]$,

$$
p_{k,n}=\frac{\lambda_n}{k b_n},
$$

\Bin where $\lambda_n\leq k_n^2 b_n$ and $\lambda_n \rightarrow \lambda$. For $0\leq s <t\leq 1$, for $n\geq 1$, we have

\begin{equation}
\Delta p_n(s,t)=\frac{\lambda_n}{b_n} \left\{\frac{1}{k_n(s)+1}+\cdots+\frac{1}{k_n(t)}\right\}
=\frac{\lambda_n (\log k_n)}{b_n} \frac{\log k_n(t) - \log k_n(s) +o(1)}{\log k_n}. \label{eqExamp1}
\end{equation}

\Bin By setting $a_n=k_n^{\varepsilon}-1$ (for $n$ large enough), by direct computations, we may check

\begin{equation}
\sup_{(s,t)\in ]0,1]^2} \left|\frac{\lambda_n\left(\log k_n(t) - \log k_n(s)\right)}{\log k_n}-\lambda(t-s) \right|\leq \log(1+a_n^{-1}), \ n\geq 1. \label{eqExamp2}
\end{equation}

\Bin The combination of Formulas \eqref{eqExamp1} and \eqref{eqExamp2} leads to (F3b).\\

\Ni \textbf{4. Non trivial examples for (F3c) holds}. \\

\Ni (a) Let $(\varepsilon_n)_{n\geq 1}$ a sequence of real numbers converging to zero and bounded by $\lambda/2>0$. Let us set, for $n\geq 1$, $k_n(t)=[k_nt]$, $k_n<n$ and

$$
p_{k,n}=\frac{\lambda + (-1)^k \varepsilon_n}{n}, \ 1\leq k \leq k_n.
$$

\Bin For any $j \in [1,k_n]$, for any $\ell \in [1, \ n-k_n]$,

$$
p_{j,n}+\cdots+p_{j+\ell,n} \in \left\{ \frac{\lambda (\ell+1) - \varepsilon_{n}}{n}, \frac{\lambda (\ell+1)}{n}, \frac{\lambda (\ell+1) + \varepsilon_{n}}{n}\right\}
$$

\Bin and hence the partial sums $p_{j+1,n}+\cdots+p_{j+\ell,n}$ are stationary.\\

\Ni (b) Let us take $k_n=n$ and

$$
b_n=\left(1 + \frac{1}{k_n+1}\right)^{-1} + \left(1 + \frac{1}{k_n+2}\right)^{-1}+\cdots + \left(1 + \frac{1}{k_n+k_n}\right)^{-1}
$$

\Bin and
$$
p_{k,n}=\frac{\lambda_n}{b_n\left(1 + \frac{1}{k+k_n}\right)}, \ 1\leq k \leq k_n,
$$

\Bin where $\lambda_n<b_n$ and $\lambda_n\rightarrow \lambda$. \\

\Ni We will show later that  $b_n=k_n(1+o(1))$. We have

$$
p_{j,n}+\cdots+p_{j+m(n)-1,n}= \frac{\lambda_n}{b_n} \times \left\{\frac{k_n+j}{(k_n+1)+j} + \cdots + \frac{k_n+j+m(n)-1}{k_n+j+m(n)}\right\}.
$$

\Bin We denote $\Delta p_n(j)=p_{j,n}+\cdots+p_{j+m(n)-1}$. Let us bound

$$
A_{j,n}=\frac{k_n+j}{(k_n+1)+j} + \cdots + \frac{k_n+j+m(n)-1}{k_n+j+m(n)}.
$$

\Bin The function $x \mapsto h(x)=(k_n+x)/(k_n+1+x)$ is non-decreasing and concave and so, the comparison of the integral $\int_{j}^{j+m(n)-1} h(x) \ dx$ and the series is as follows:

$$
A_{j,n} - \frac{k_n+j+m(n)-1}{k_n+j+m(n)}  \leq \int_{j}^{j+m(n)-1} \frac{k_n+x}{k_n+1+x} \ dx\leq  A_{j,n} - \frac{k_n+j}{k_n+j+1}.
$$

\Bin But

$$
\int_{j}^{j+m(n)-1} \frac{k_n+1+x}{k_n+x} \ dx=\int_{j}^{j+m(n)-1} \left(1 + \frac{1}{k_n+x}\right) \ dx
$$

\Bin and by using the primitive

$$
\left(x + \log (k_n+x)\right) \ of \ \left(1 + \frac{1}{k_n+x}\right),
$$

\Bin we get

\begin{eqnarray}
A_{j,n} - \frac{k_n+j+m(n)-1}{k_n+j+m(n)} &\leq& (m(n)-1) + \log \frac{k_n+j+m(n)-1}{k_n+j} \notag\textbf{}\\
&\leq&  A_{j,n} - \frac{k_n+j}{k_n+j+1}. \label{bound10}
\end{eqnarray}

\Bin But all the terms, except $A_{j,n}$, when divided by $b_n \sim k_n$, go to zero uniformly in $j$. Let show this for one them, the middle term for example,

$$
\frac{1}{k_n} \log \frac{k_n-1}{2k_n} \leq B_{j,n}=\frac{1}{k_n} \log \frac{k_n+j+m(n)-1}{k_n+j} \leq \frac{1}{k_n} \log \frac{3k_n-1}{k_n}
$$

\Bin which shows that $\sup_{1\leq j \leq k_n} B_{j,n} \rightarrow 0$. By using the same remarks, ideas and methods, we have

$$
b_n-\frac{2 k_n}{2k_n+1} \leq \int_{1}^{k_n} \left(1 - \frac{1}{k_n+x+1}\right) \ dx \leq b_n-\frac{k_n+1}{k_n+2}.
$$

\Bin The integral is $(k_n-1) - \log( (2k_n+1)/(k_n+2))$ and so $b_n=k_n (1+o(1))$. Now, Formula \eqref{bound10} becomes

$$
\frac{A_{j,n}}{b_n}+o(1) \leq \frac{m(n)}{k_n} (1+o(1)) + o(1)\leq \frac{A_{j,n}}{b_n}+o(1),
$$

\Bin where all the small $o's$ are uniform in $j$. We already know that

$$
\left|\frac{m(n)}{k_n}-(t-s)\right|\leq \frac{1}{k_n}.
$$

\Bin The combination of all these facts gives

$$
\sup_{1\leq j \leq k_n} \left|\Delta p_n(j) - \lambda (t-s)\right| \rightarrow 0
$$

\Bin and so Condition (F3c) holds. $\blacksquare$

\section{Conclusion} 
\Ni We saw how a very simple result in probability that can be proved by moment generating function can become sophisticated for sums of independent but non-stationary random variables. In this paper, we treated the question in the general frame of the weak laws in $\ell^{\infty}(0,1)$ and learned how to deal with it in a general frame to have new results. The most important thing is the construction of the frame that will allow extensions for non-stationary and dependent data, for associated data in the first place.\\

\Ni \textbf{Acknowledgment}. The authors wish to thank the members of Imhotep Mathematical Center (IMC) for their comments.

\end{document}